\documentclass[12pt]{article}
\usepackage{amsmath}
\usepackage{amsfonts}
\usepackage{amssymb}
\usepackage{epsfig}
\newtheorem{lem}{Lemma}

\newtheorem{thm*}{Theorem}
\newtheorem{cor}{Corollary}

\newtheorem{rem}{Remark}
\newtheorem{pr}{Proposition}

\renewcommand{\a}{\alpha}
\renewcommand{\b}{\beta}
\newcommand{\g}{\gamma}
\renewcommand{\l}{\lambda}
\newcommand{\f}{\varphi}
\newcommand{\ind}{{\rm ind}}

\newcommand{\HH}{\mathbb{H}}
\newcommand{\OO}{\mathbb{O}}
\newcommand{\Z}{\mit{Z}}
\newcommand{\Cay}{{\rm Cay\,}}
\newcommand{\Ca}{{\rm Cay}}
\newcommand{\Reg}{{\rm Reg\,}}

\def\bes{\begin{eqnarray*}}
\def\ees{\end{eqnarray*}}
\def\bee{\begin{eqnarray}}
\def\eee{\end{eqnarray}}

\newcommand{\G}{{\mit\Gamma}}

\newcommand{\la}{\langle}
\newcommand{\ra}{\rangle}
\begin{document}

\centerline{\Large{{\bf Alternative $\mit{M_2}$-algebras and $\G$-algebras}}} 

\vspace{7mm}

\centerline{\bf\large A.\,Grishkov,  I.\,Shestakov}

\vspace{6mm}

\begin{abstract}
Recently V.\,H.\,L\'opez Sol\'is and I.\,Shestakov \cite{LSSh} solved an old problem by N.Jacobson \cite{J1} on describing of unital alternative algebras containing the $2\times 2$ matrix algebra $M_2$ as a unital subalgebra.
Here we give another description of $M_2$-algebras via the 6-dimensional alternative superalgebra $B(4,2)$ and an auxiliar $\mathbf {\mit Z}_2$-graded algebra $\mit\Gamma$. It occurs that the category of alternative $M_2$-algebras is isomorphic to the category of $\G$-algebras.  For any associative and commutative algebra $A$, we give a construction of a $\G$-algebra $\G(A)$,  which turns to be  a Jordan superalgebra; if $A$ is a domain then $\G(A)$ is a prime superalgebra.  We describe also the free $\G$-algebras and construct their bases.
\end{abstract}
\begin{flushright}
\sl{To the memory of our dear friend Sasha Anan'in}
\end{flushright}
\section{Introduction.}

The classical Wedderburn Coordinatization Theorem says that if a unital associative algebra $A$ contains a matrix algebra $M_n(F)$  over a field $F$ with the same identity element then it is itself a matrix algebra, $A\cong M_n(D)$, ``coordinated'' by $D$. Generalizations and analogues of this theorem were proved for various classes of algebras and superalgebras \cite{J1,K,LDS1,LDS2,  MSZ,CE,Mac,PS,Sc}. The common content of all these results is that if an algebra (or superalgebra) contains a certain subalgebra (matrix algebra, octonions, Albert algebra) with the same unit then the algebra itself has the same structure,  but not over the basic field rather over a certain algebra that ``coordinatizes'' it.  The Coordinatization Theorems play important role in structure theories, especially in classification theorems, and also in the representation theory, since quite often an algebra $A$ coordinated by $D$  is Morita equivalent to $D$, though they could belong to different classes (for instance, Jordan algebras are coordinated by associative and alternative algebras).

I.\,Kaplansky \cite{K} proved an analogue  of Wedderburn's theorem for alternative algebras containing the Cayley algebra. He showed
that if $A$ is an alternative algebra with identity element $1$ which contains a subalgebra ${B}$
isomorphic to a Cayley algebra and if $1$ is contained in ${B}$, then $A$ is isomorphic to the Kronecker product ${B}\otimes T$, where $T$ is the center of $A$.

The Wedderburn coordinatization theorem in the case $n\geq 3$ admits a generalization for  alternative algebras, since every alternative algebra $A$ which contains a subalgebra $M_n(F)$ $(n\geq 3)$ with the same identity element is associative (see \cite[Corollary $11$, Chapter $2$]{Sc}). The result is not true for $n=2$, the split Cayley algebra and its 6-dimensional subalgebra are counterexamples. The problem of description of alternative algebras containing $M_2(F)$ or, more generally,  a generalized quaternion algebra $\mathbb{H}$ with the same identity element was posed by Jacobson \cite{J1}.

\smallskip

 In \cite{LSSh},  this problem was solved for the split case $\mathbb{H}\cong M_2(F)$. The corresponding  $M_2(F)$-coordinatization  in \cite{LSSh}  involves two ingredients:  an associative algebra $D$ and  a commutative $D$-bimodule $V$ (that is, $V$ is annihilated by any commutator of elements of $D$), on which a skew-symmetric mapping is defined with values in the center of $D$, satisfying Pl\"ucker relations. More exactly,  $A=M_2(D)\oplus V^2$, with a properly defined multiplication. 

\smallskip

Here we give another characterization of $M_2(F)$-algebras, based on the 6-dimensional simple alternative superalgebra $B=B(4,2)$ \cite{Sh1} and an auxiliar $\mathbf {\mit Z}_2$-graded algebra $\G$: a unital  alternative algebra  $A$ is an  $M_2(F)$-algebra if and only if $A$ is a {\em$\G$-envelope} of the superalgebra $B$:  $A=\G_0\otimes B_0+\G_1\otimes B_1$.  Moreover, the category of alternative $M_2$-algebras is isomorphic to the category of $\G$-algebras.  For any associative and commutative algebra $A$, we give a construction of a $\G$-algebra $\G(A)$,  which turns to be  a Jordan superalgebra; if $A$ is a domain then $\G(A)$ is a prime superalgebra.  Finally,  we describe  free $\G$-algebras and construct their bases.  It occurs that these algebras are closely related to coordinate algebras of grassmannians $Gr(2,n)$.  

\smallskip
Throughout this paper the ground field F is of arbitrary characteristic.

\section{Definitions, examples, and preliminary results}

Let $A$ be a composition algebra (see \cite{J2, KS, Sc, ZSSS}). Recall that $A$ is a unital alternative algebra, it has an involution $a\mapsto \bar a$ such that the {\em {trace}} $t(a)=a+\bar a$ and {\em {norm}} $n(a)=a\bar a$ lie in $F$.

An alternative bimodule $V$ over a composition algebra $A$ is called a \textit{{Cayley bimodule}} if it satisfies the relation
\begin{equation}\label{e3}
av=v{\bar a},
\end{equation}
where $a\in A$, $v\in V$, and $a\rightarrow {\bar a}$ is the canonical involution in $A$.

\smallskip
Typical examples of composition algebras are the algebras of (generalized) quaternions $\HH$ and octonions $\OO$ (or {\em a Cayley algebra}) with symplectic involutions.
 Recall that $\OO=\HH\oplus v\HH$, with the product defined by
 \bee\label{CD_product}
{{a\cdot b=ab,\ a\cdot vb=v(\bar ab),\ vb\cdot a=v(ab), \ va\cdot vb= (b\bar a)v^2,}}
\eee
where $a,b\in\HH,\ 0\neq v^2\in F, \ a\mapsto \bar a$ is the symplectic involution in $\HH$.

The subspace {{$v\HH\subset\OO$}} is invariant under multiplication by elements of $\HH$ and it gives an example of a {{ Cayley bimodule}} over $\HH$.
If  {{$\HH$ is a division algebra then $v\HH$ is irreducible, otherwise $\HH\cong M_2(F)$ and  
\[
v\HH=\la ve_{{22}},-ve_{{12}}\ra\oplus\la -ve_{{21}},ve_{11}\ra,
 \] 
 where $M_2(F)$-bimodules $\la ve_{{22}},\,-ve_{{12}}\ra$ and $\la -ve_{{21}},ve_{11}\ra$ are both isomorphic to the 2-dimensional Cayley bimodule
$\Ca=F\cdot m_{1}+F\cdot m_{2}$, with the action of $M_{2}(F)$ given by
\bee\label{id_cay}
e_{ij}\cdot m_{k}=\delta_{ik}m_{j},\ \ \ \  
m\cdot a={\bar a}\cdot m, 
\eee
where $a\in M_{2}(F),\, m\in \Ca,\ i,j,k\in\{1,2\}$ and $a\mapsto {\bar a}$ is the symplectic involution in $M_{2}(F)$. In the last case the algebra $\OO=M_2(F)\oplus vM_2(F)$ is called {\em the split octonion algebra}.

 \smallskip
 Let us call a unital alternative algebra $A$ an {{\it$M_2$-algebra}} if there exists a homomorphism of unital algebras $\phi:M_2\rightarrow A$, where $M_2=M_2(F)$ .
\vspace{5mm}
 
{\bf{Examples of $M_2$-algebras:}}
\begin{enumerate}
\item  $A$ associative, $A\supseteq M_2\ni1_A\Rightarrow A=M_2(B)\cong M_2\otimes B,\ B$ associative.
\item $\OO$ split octonion algebra, {$ \OO=M_2\oplus (M_2)v$.}
  \item {$S=M_2(F)\oplus \Ca\subseteq\OO,\  \Ca^2=0$,} the split null extension of $M_2$ by bimodule $\Ca$. 
  \item {$G(B(4,2))=G_0\otimes M_2+G_1\otimes \Ca,\ char\,F=3$},
  the Grassmann envelope of the simple alternative superalgebra ${B(4,2)=M_2\oplus \Ca}$ (see \cite{Sh1}), with the following  multiplication in ${\Cay}$:
  \bes
  {m_1^2 = e_{21},\  m_2^2 = -e_{12}, \ m_1\,m_2= -e_{11},\ m_2\,m_1= e_{22}.}
\ees
  \end{enumerate}
 \begin{rem}
 The odd product in $B(4,2)$ in \cite{Sh1} has different sign; one can get the old product by the following change of the basis: 
 \bes
 e_{21}\leftrightarrow -e_{21}.\ e_{12}\leftrightarrow -e_{12}, \ m_{2}\leftrightarrow -m_{2}.
 \ees 
 \end{rem}

Any $M_2$-algebra $A$ may be considered as a unital alternative $M_2$-bimodule. The structure of such bimodules is given by the following result:

\begin{thm*}\label{thm1}
\cite{J2, Sh1}
 Let $V$ be a unital alternative $M_2$-bimodule. Then $V$ is completely reducible, moreover, $V=V_a\oplus V_c$, where $V_a$ is an associative $M_2$-bimodule and is a direct some of regular bimodules $\Reg M_2$, while $V_c$ is a Cayley $M_2$-bimodule which is a direct sum of irreducible Cayley bimodules of type  $\Cay.$
\end{thm*}


In particular, any $M_2$-algebra $A$ may be written as ${A=A_a\oplus A_c}$. It was proved in \cite{LSSh} that 
\bes 
{  {A_cA_a+A_aA_c \subseteq A_c,\ A_aA_a\subseteq A_a, \ A_cA_c\subseteq A_a.}}
\ees
That is, $A$ is a $\Z_2$-graded algebra. 
We have {$A_a=(\Reg\,M_2)^n,\ A_c=(\Cay)^m$}. 
Moreover, the subalgebra $A_a$ is associative  \cite{LSSh}.


The ${M_2(F)}$-coordinatization  in \cite{LSSh} involves two ingredients: an alternative ${M_2(F)}$-algebra $A$
is ``coordinated'' by an associative algebra $D$ and by a commutative $D$-bimodule ${V}$ (that
is, ${V}$ is annihilated by any commutator of elements of ${D}$), on which a skew-symmetric form
 is defined with values in the center of ${D}$, satisfying Pl\"ucker relations. 
More exactly,
\[{A = M_2(D)\oplus V^2,}\] with a properly defined multiplication.

In more details,  let ${D}$ be an associative unital algebra and ${V}$ be a left ${D}$-module such that ${[D,D]}$ annihilates ${V}$. Clearly, in this case ${V}$ has a structure of a commutative ${D}$-bimodule with ${v\cdot a=a\cdot v,\ v\in V,\,a\in D}$. Assume that there exists a ${D}$-bilinear skew-symmetric mapping ${\la,\ra :V^2\rightarrow D}$ such that ${\la V,V\ra\subseteq Z(D)}$ and for any ${u,v,w\in V}$
\bee
{\la u,v\ra w+\la v,w\ra u+\la w,v\ra u=0.\label{1}}
\eee

Consider ${A=M_2(D)\oplus V^{2}}$. 
Let {$X,Y\in A,$ $X=X_a+(x,y),\ Y=Y_a+(z,t)$}, where {$X_a,\ Y_a\in M_2(D)$} and {$(x,y), (z,t)\in V^2$}. Define  a product in $A$ by formula:
\bes
{XY=X_aY_a+\left(\begin{array}{cc}
-\la x,t\ra&-\la y,t\ra\\
\la x,z\ra&\la y,z\ra
\end{array}\right)+(z,t)X_a+(x,y)(Y_a)^*,}
\ees
where {$\left(\begin{array}{cc}
a&b\\
c&d
\end{array}\right)^*=\left(\begin{array}{cc}
d&-b\\
-c&a
\end{array}\right)$.}

\begin{thm*}\label{thm2} \cite{LSSh}  
{The algebra $A$ with the product defined above is an alternative unital algebra containing $M_2(F)$ with the same identity element. Con\-ver\-se\-ly, eve\-ry uni\-tal alter\-nati\-ve algebra that con\-tains the mat\-rix al\-geb\-ra $M_2(F)$ with the same iden\-tity element has this form.}
 \end{thm*}
 
 \section{The new construction}

Consider the vector space direct sum ${\G}=D\oplus V$, where $D$ and $V$ are taken from theorem 2,  and define a multiplication on it as follows
\bes
(a+u)(b+v)=(ab+\la u,v\ra)+(av+bu).
\ees
Then ${\G}$ becomes a ${\Z}_2$-graded algebra with ${\G}_{0}=D,\, {\G}_1=V$, that satisfied the following conditions:
\begin{itemize}
 \item[(i)]  ${\G}_0$ is a unital associative algebra, $[{\G}_0,{\G}_1]=(\G_{0},\G,\G)=0$,
\item[(ii)]   ${\G}_1^2\subseteq Z({\G}_0)$,
\item[(iii)] $x y+yx=0,\ x,y\in{\G_1}$, 
\item[(iv)] $(xy)z+(yz)x+(zx)y=0,\ x,y,z\in{\G_1}$.
\end{itemize}

\smallskip

\begin{thm*}\label{thm3}
The ``${\G}$-envelope'' ${\G}(B(4,2))={\G}_0\otimes M_2+{\G}_1\otimes Cay$ is isomorphic to the algebra $M_2({\G}_0)\oplus {\G}_1^2$ from theorem 2. In particular, an algebra $A$ is an alternative $M_2$-algebra if and only if $A={\G}(B(4,2))$ for a certain $\Z_2$-graded algebra ${\G}$ satisfying the above conditions.
\end{thm*}
\begin{rem}
Note that the superalgebra $B(4,2)$ is alternative only in $char\,F= 3$ case, but the theorem holds in any characteristic.
\end{rem}
{\em\bf Proof of the theorem.}
It is clear that ${\G}_0\otimes M_2\cong M_2({\G}_0)$. Let us prove that the mapping
 \bes
 \varphi: X+(x,y)\mapsto X+x\otimes m_1+y\otimes m_2,\ X\in M_2({\G}_0),\, x,y\in {\G}_1,
 \ees
 is an isomorphism of $M_2({\G}_0)\oplus {\G}_1^2$ and ${\G}(B(4,2))$. 
Let $X=\left(\begin{array}{cc}
d_{11}&d_{12}\\
d_{21}&d_{22}
\end{array}\right), d_{ij}\in {\G}_1$. Consider
\bes
\varphi (X\cdot (x,y))&=&\varphi (xd_{11}+yd_{21},xd_{12}+yd_{22})\\
&=&(xd_{11}+yd_{21})\otimes m_1+(xd_{12}+yd_{22})\otimes m_2.
\ees
On the other hand,
\bes
\varphi (X)\varphi (x,y)&=&X(x\otimes m_1+y\otimes m_2)=(\sum_{ij}d_{ij}\otimes e_{ij})(x\otimes m_1+y\otimes n)\\
&=&(d_{11}x+d_{21}y)\otimes m_1+(d_{12}x+d_{22}y)\otimes m_2.
\ees
Since $[{\G}_0,{\G}_1]=0$, we have  $\varphi (X\cdot (x,y))=\varphi (X)\varphi (x,y)$. Similarly, $\varphi ( (x,y)\cdot X)=\varphi (x,y)\varphi (X)$.
Let now $z,t\in {\G}_1$, consider 
\bes
\varphi ((x,y)(z,t))&=&\varphi \left(\begin{array}{cc}
-\la x,t\ra&-\la y,t\ra\\
\la x,z\ra&\la y,z\ra
\end{array}\right)=\left(\begin{array}{cc}
-\la x,t\ra&-\la y,t\ra\\
\la x,z\ra&\la y,z\ra
\end{array}\right).
\ees
On the other hand,
\bes
\varphi (x,y)\varphi (z,t)&=&(x\otimes m_1+y\otimes n)(z\otimes m_1+t\otimes m_2)\\
&=&\la x,z\ra \otimes e_{21}-\la x,t\ra \otimes e_{11}+\la y,z\ra e_{22}-\la y,t\ra\otimes e_{12}\\
&=&\left(\begin{array}{cc}
-\la x,t\ra&-\la y,t\ra\\
\la x,z\ra&\la y,z\ra
\end{array}\right).
\ees
This proves the theorem.

\hfill$\Box$

Let us consider the examples of $M_2$-algebras from section 2 and determine the structure of the corresponding algebra ${\G}$ in every case.
\begin{enumerate}
\item  $A$ associative, $A\supseteq M_2\ni1_A\Rightarrow A=M_2(B)\cong M_2\otimes B,\ B$ associative.
In this case ${\G}={\G}_0=B$.
\item $\OO$ split octonion algebra, {$ \OO=M_2\oplus (M_2)v$.} Here ${\G}=B(1,2)$, the 3-dimensional simple superalgebra from \cite{Sh1} which is alternative in characteristic 3 case. $B(1,2)_0=F,\ B(1,2)_1=Fx+Fy,\ x^2=y^2=0, xy=-yx=1$.
  \item {$S=M_2(F)\oplus \Ca\subseteq\OO,\  \Ca^2=0$,} the split null extension of $M_2$ by bimodule $\Ca$. In this case ${\G}=F+Fx,\, x^2=0,$ is the 2-dimensional algebra with ${\G}_0=F,\, {\G}_1=Fx$.
  \item {$G(B(4,2))=G_0\otimes M_2+G_1\otimes \Cay,\ char\,F=3$},
  the Grassmann envelope of the simple alternative superalgebra $B(4,2)$. Here ${\G}=G$ since the Grassmann algebra  $G$ satisfies the conditions for ${\G}$ in the case of characteristic 3.
  \end{enumerate}

\section{Tensor algebras of bimodules and free ${\G}$-algebras}

Recall the definition of a tensor algebra of bimodule (see \cite{KOS}). Let $A$ be an algebra in a variety $\cal M$ and let $V$ be an $\cal M$-bimodule over $A$. Consider the free $\cal M$-algebra $F_{\cal M}[A\oplus V]$ and let $I$ be the ideal of this algebra generated by the set $\{a*b-ab,\, a*v-a\cdot v,\, v*a-v\cdot a\,|\, a,b\in A,\, v\in V\}$, where $*$ denotes the multiplication in 
$F_{\cal M}[A\oplus V]$ and $a\cdot v, \,v\cdot a$ denote the action of $A$ on $V$. Then the quotient algebra $F_{\cal M}[A\oplus V]/I$ is called {\em the tensor algebra of the $A$-bimodule $V$}.

By the standard arguments, one can  prove the following universal property of tensor algebra.
\begin{pr}
Let $B\in\cal M$ and let $\f:A\rightarrow B$ be a homomorphism of algebras. Then $B$ has a natural structure of an $A$-bimodule. Now, for any homomorphism of $A$-bimodules $\psi:V\rightarrow B$ there exists a unique homomorphism of algebras 
$\tilde\psi : F_{\cal M}[A\oplus V]\rightarrow B$ such that $\tilde\psi(a)=\f(a),\, \tilde\psi(v)=\psi(v)$ for any $a\in A,\,v\in V$.
\end{pr}

In particular, the tensor algebra $M_2[V]$ of an alternative $M_2$-bimodule $V$ plays a role of a free object in the category of alternative $M_2$-algebras: for any $M_2$-algebra $B$, any homomorphism of $M_2$-bimodules $\f:V\rightarrow B$ is uniquely extended to an algebra homomorphism $\tilde\f:M_2[V]\rightarrow B$.
 
 \smallskip
 Let us call a $\Z_{2}$-graded algebra satisfying conditions (i) - (iv) of section 3 a {\em ${\G}$-algebra}.
 We want to prove 
 \begin{thm*}\label{thm4} 
 The category of alternative $M_{2}$-algebras is isomorphic to the category of ${\G}$-algebras.
 \end{thm*}
 {\em \bf Proof}. 
 There is a natural functor from the category of $\G$-algebras to the category of $M_{2}$-algebras:
 $F:\G\rightarrow  \G(B(4,2))$, which sends a morphism of $\G$-algebras $\f:\G\rightarrow \G'$ to the morphism $F(\f):\G(B(4,2))\rightarrow  \G'(B(4,2))$ identical on $B(4,2)$. 
 
 It is clear from the proof of theorem \ref{thm3} that this functor is bijective on objects: every $M_{2}$-algebra $A=M_{2}(D)\oplus V^{2}$ defines uniquely $\G_0=D$ and $\G_{1}=V$.  It remains to show that any morphism $\f:\G(B(4,2))\rightarrow  \G'(B(4,2))$ is induced by a morphism $\psi:\G\rightarrow \G'$ such that $\f=F(\psi)$. Denote $A=\G(B(4,2)),\ A'=\G'(B(4,2))$. Since $\f$ is identical on $M_{2}$, it is a homomorphism of $M_{2}$-bimodules; in particular, $\f(A_{a})=\f(M_{2}(\G_{0}))\subset A'_{a}=M_{2}(\G'_{0})$ and $\f(A_{c})=\f(\G_{1}\otimes \Ca)\subset (A')_{c}=\G'_{1}\otimes \Ca$.  It is well known that a homomorphism of matrix algebras is induced by a homomorphism of their coordinates,  hence there exists a homomorphism $\psi_{0}:\G_{0}\rightarrow \G'_{0}$ which induces $\f|_{A_{a}}$.
 
 Now, fix $0\neq\g\in\G_{1}$. Let $\f(\g\otimes m_{1})=\a\otimes m_{1}+\b\otimes m_{2}, \ \f(\g\otimes m_{2})=\l\otimes m_{1}+\mu\otimes m_{2}$ for some $\a,\b,\l,\mu\in \G'_{1}$. We have 
 \bes
\a\otimes m_{1}+\b\otimes m_{2}&=& \f(\g\otimes m_{1})=\f((1\otimes e_{11})(\g\otimes m_{1}))\\
&=&(1\otimes e_{11})\f(\g\otimes m_{1})=\a\otimes m_{1},
 \ees
 which implies $\b=0$. Similarly, $\l=0$. Futhermore,
 \bes
 \mu\otimes m_{2}&=&\f(\g\otimes m_{2})=\f((1\otimes e_{12})(\g\otimes m_{1}))\\
 &=&(1\otimes e_{12})\f(\g\otimes m_{1})=\a\otimes m_{2},
 \ees
 which implies $\a=\mu$. Therefore, $\f(\g\otimes\Ca)=\a\otimes\Ca$, and we put $\psi_{1}(\g)=\a$. 
 One can easily check that $\psi=\psi_{0}+\psi_{1}:\G\rightarrow\G'$ is a homomorphism of algebras such that $F(\psi)=\f$.
 
 \hfill$\Box$
 
 \begin{cor}\label{cor1}
 Let ${\G}[X_{0},X_{1}]$ be the free $\G$-algebra on  sets $X_{0}$ and $X_{1}$ of even and odd generators. Then ${\G}[X_{0},X_{1}](B(4,2))\cong M_{2}[V]$, where $V=Reg^{\#X_{0}}\oplus \Ca^{\#X_{1}}.$
\end{cor}

In view of the Corollary, it seems important to determine the structure of  free $\G$-algebras.
 
 \smallskip
 
 Let $V=F^n$ be an $n$-dimensional vector space over $F$ and $F[Gr_{2}(V)]$ be the coordinate algebra of the Grassmannian $Gr(2,n)=Gr(2,V)$.
 Recall that $F[Gr(2,V)]\cong F[V^{\wedge 2}]/P$, where $P$ is the ideal generated by the {\em double Pl\"ucker relations}
 \bes
 (u\wedge v)(w\wedge z)+(u\wedge w)(z\wedge v)+(u\wedge z)(v\wedge w),\ u,v,w,z\in V.
 \ees
 Furthermore, consider the tensor product $ F[Gr(2,V)]\otimes V$ which has a natural structure of $F[Gr(2,V)]$-module.
 Denote by $F[Gr(2,V)]_1$  the quotient module  $(F[Gr(2,V)]\otimes V)/I$ where $I$ is the  $F[Gr(2,V)]$-submodule generated by the {\em ordirnary  Pl\"ucker relations}
  \bes
 (u\wedge v)\otimes w+(v\wedge w)\otimes u+(w\wedge u)\otimes v,\ u,v,w\in V.
 \ees
Define  multiplication in $F[Gr(2,V)]_1$ with results in $F[Gr(2,V)]$ by setting
\bes
(a\otimes u)(b\otimes v)=ab(u\wedge v)\otimes 1,
\ees
where $a,b\in F[Gr(2,V)],\ u,v\in V$. The product is defined correctly, since
\bes
 ((u\wedge v)\otimes w+(v\wedge w)\otimes u+(w\wedge u)\otimes v)(a\otimes z)&=&\\
a( (u\wedge v)(w\wedge z)+(v\wedge w)(u\wedge z)+(w\wedge u)(v\wedge z))\otimes 1&=&0.
\ees
\begin{thm*}\label{thm5}
The $\Z_2$-graded algebra $\G[\emptyset;V]=F[Gr(2,V)]+F[Gr(2,V)]_1$ is a free $\G$-algebra generated by the space of odd generators $V$.
\end{thm*}
{\em \bf Proof}. In fact, it is easy to check that the unital algebra $\G[\emptyset;V]$ satisfies conditions (i) -- (iv) defining  $\G$-algebras, and is generated by the space $V$.  Let $\G=\G_0+\G_1$ be a $\G$-algebra and $\f:V\rightarrow \G_1$ a linear mapping. For $v\in V$ denote by $\bar v$ its image in $\G_1$. The mapping $\f$ is extended to a linear mapping $V^{\wedge 2}\rightarrow \G_0,\ u\wedge v\mapsto \bar u\bar v$ and further to an algebra homomorphism $F[V^{\wedge 2}]\rightarrow \G_0$.  In view of condition (iv) the ordinary Pl\"ucker relation holds in $\G$. Moreover, for any $u,v,w,z\in V$ we have in $\G$
\bes
(\bar u\bar v)(\bar w\bar z)+(\bar u\bar w)(\bar z\bar v)+(\bar u\bar z)(\bar v\bar w)=
\bar u(\bar v\bar w\bar z+\bar w\bar z\bar v)+\bar z\bar v\bar w)=0.
\ees
Therefore, the mapping $\f$ can be extended to an algebra homomorphism $\tilde\f:\G[\emptyset;V]\rightarrow \G$.

\hfill$\Box$

Let now $U$ be another vector space, construct the free $\G$-algebra $\G[U;V]$ generated by the space $U$ of even generators and the space $V$ of odd generators.  Denote by $F\la U\ra$ and by $F[U]$ the free associative and polynomial algebras over the space $U$. Furthermore, by $F[Gr(2,V)]^{0}$ denote the augmentation ideal of the algebra $F[Gr(2,V)]$, that is, the ideal of elements without scalar terms.
Consider the $\Z_2$-graded vector space 
\bes
\G[U;V]=(F\la U\ra \oplus (F[U]\otimes F[Gr(2,V)]^{0}))\oplus (F[U]\otimes F[Gr(2,V)]_1),
\ees
with $\G[U;V]_0=F\la U\ra \oplus (F[U]\otimes F[Gr(2,V)]^{0})$ and $\G[U;V]_1=F[U]\otimes F[Gr(2,V)]_1$.
Observe that 
\bes
I= (F[U]\otimes F[Gr(2,V)]^{0})\oplus (F[U]\otimes F[Gr(2,V)]_1)=\\
 F[U]\otimes (F[Gr(2,V)]^{0}+ F[Gr(2,V)]_1)\subseteq F[U]\otimes \G[\emptyset;V]
 \ees
Define multiplication on $\G[U;V]$ in the following way: the space  
 $I$ is an ideal of $\G(U;V)$ with the product defined as in a subalgebra of the algebra 
$F[U]\otimes \G[\emptyset;V]$;  the algebra $F\la U\ra$ is a subalgebra of $\G[U;V]$, and the element $f\in F\la U\ra$ acts on $I$ by
\bes
f\cdot (g\otimes a+h\otimes\ b\otimes v)=\bar f g\otimes a+\bar f h\otimes\ b\otimes v,
\ees
where $g,h\in F[U],\,a,b\in F[Gr(2,V)],\, v\in V,$ and $\bar f\in F[U]$ is the image of $f$ under the natural epimorphism $F\la U\ra\rightarrow F[U]$.
\begin{thm*}\label{thm6}
The algebra $\G[U;V]$ with the multiplication defined above is a free $\G$-algebra generated by the spaces $U$ and $V$ of even and odd generators.
\end{thm*}
{\em \bf Proof}. First of all, one can easily check that the algebra $\G[U;V]$ satisfies  conditions (i) -- (iv). Furthermore, let $\G=\G_0+\G_1$ be a $\G$-algebra and $\f:U\rightarrow \G_0,\ \psi:V\rightarrow \G_1$ be linear mappings. By above,  $\psi$ can be extended to  an algebra homomorphism $\tilde\psi:\G[\emptyset;V]\rightarrow \G$.   By the property of free algebras, there exists also an algebra homomorphism $\tilde\f :F\la U\ra\rightarrow \G_0$ extending $\f$. 
Now the mapping 
\bes
(\tilde\f+\tilde\psi):f+(h\otimes a+g\otimes b\otimes v)\mapsto \tilde\f(f)+\tilde\f(h)\tilde\psi(a)+\tilde\f(g) \tilde\psi(b)\psi(v)
\ees
 for $f\in F\la U\ra,\ g,h\in F[U],\, a,b\in F[Gr(2,V)]$ and $v\in V$ is an algebra homomorphism of $\G[U;V]$ to $\G$ extending $\f+\psi$.

\hfill$\Box$

\section{$\G$-algebras and Jordan superalgebras}

Observe that if $\G$ is a $\G$-algebra with commutative even part $\G_0$ then $\G$ is a commutative superalgebra. Moreover, in this case it is a Jordan superalgebra.

\begin{pr}\label{pr2}
Let $\G$ is a $\G$-algebra with commutative even part $\G_0$. Then $\G$ is a Jordan superalgebra.  Moreover, if $char\, F=3$ then $\G$ is an alternative superalgebra.
\end{pr}
{\bf Proof.} Recall that a commutative superalgebra is called a Jordan superalgebra if it satisfies the super-identity
\bee\label{SJ}
(xy,z,t)+(-1)^{\bar y\bar z+\bar y\bar t+\bar z\bar t}(xt,z,y)+(-1)^{\bar x(\bar y+\bar z+\bar t)+\bar z\bar t}(yt,z,x)=0,
\eee
where $\bar x$ for $x\in \G_0\cup\G_1$ denotes the parity of element  $x$: $\bar x=i \Leftrightarrow x\in\G_i$.
Note that $[\G_0,\G]=0$, hence $\G_0$ is contained in the center of $\G$. Therefore, if at least 2 elements of $x,y,z,t$ lie in $\G_0$ or $z\in \G_0$, all the associators in identity \eqref{SJ} vanish. If all the elements $x,y,z,t$ are in $\G_1$ then $xy,\, xt,\, yt\in\G_0$, and again  \eqref{SJ} holds. Therefore, it suffices to consider the case when  $x,y,z\in\G_1,\, t=a\in \G_0$. We have
\bes
&&(xy,z,a)-(xa,z,y)+(ya,z,x)=-(xa,z,y)+(ya,z,x)\\
&=&a(-(x,z,y)+(y,z,x))=a(-xz\cdot y+x\cdot zy+yz\cdot x-y\cdot zx)\\
&=&a(-(xz+zx)a+x(zy+yz))=0.
\ees

Furthermore, let $char\,F=3$. Since $\G_0\subseteq Z(\G)$, in order to check alternatively we have to consider only associators on odd generators.
Let $x,y,z\in\G_1$, then we have
\bes
(x,y,z)-(x,z,y)&=&xy\cdot z-x\cdot yz-xz\cdot y+x\cdot zy=xy\cdot z-2yz\cdot x+zx\cdot y\\
&=&xy\cdot z+yz\cdot x+zx\cdot y=0,
\ees
and similarly $(x,y,z)-(y,x,z)=0$, hence the superalgebra $\G$ is alternative.

\hfill$\Box$

An important example of supercommutative $\G$-algebras  can be obtained as follows. Let $A$ be a commutative associative algebra, consider
$\G(A)=A\oplus A^2$ with the grading $\G(A)_0=A,\, \G(A)_1=A^2$ and the following multiplication:
\bes
a\cdot b=ab,\, a\cdot (b,c)=(b,c)\cdot a=(ab,ac),\, (a,b)\cdot (c,d)=ad-bc;\ a,b,c,d\in A.
\ees
We have only to check condition (iv) in the definition of $\G$-algebra. Consider
\bes
&(a,b)(c,d)\cdot (e,f)+(c,d)(e,f)\cdot (a,b)+(e,f)(a,b)\cdot (c,d)&\\
&=(ad-bc)(e,f)+(cf-de)(a,b)+(eb-fa)(c,d)=(0,0).&
\ees 
Therefore, $\G(A)$ is a $\G$-algebra. Since $A$ is commutative, $\G(A)$ is a Jordan superalgebra.

\begin{pr}\label{pr3}
Let $A$ be a domain, then  $\G(A)$ is a central order in the simple Jordan superalgebra of type $B(1,2)$ (the superalgebra of a skew-symmetric bilinear form on a 2-dimensional vector space). In particular, in this case $\G(A)$ is prime and special.
\end{pr}
{\bf Proof.} In fact, let $K$ be the quotient field of $A$, then we have an inclusion $\G(A)\subseteq\G(K)$; moreover, since $K=A^{-1}A$ and $A=Z(\G(A) )$, we have $A^{-1}\G(A)=\G(K)$.  It is clear that $\dim_K\G(K)=3$ and $\G(K)=B(1,2)$ as a $K$-superalgebra. 

It is well known that $B(1,2)$ is a special superalgebra, hence so is $\G(A)$. Finally, a central order in a simple (super)algebra is evidently  prime.

\hfill$\Box$

\section{Bases of free $\G$-algebras}

In this section we will construct bases of free $\G$-algebras defined in terms of free generators.
 
 Let $A_{n}=F[x_{1},\ldots,x_{n};y_{1},\ldots,y_{n}],$ consider the $\G$-algebra $\G(A_n)$.
 We want to prove that the subalgebra of $\G(A_n)$ generated by the odd elements $v_i=(x_i,y_i),\, i=1,\ldots,n,$ is a free $\G$-algebra on these set of generators.
 
 Denote $\a_{ij}=v_iv_j=x_{i}y_{j}-x_{j}y_{i}\, (1\leq i<j\leq n), 
 \ S_{n}=F[\a_{12},\ldots,\a_{(n-1)n}]\subset A_{n},\ \ V_i=Fv_i,\, V=\sum_{i=1}^nV_i$. 
 We have the relations
 \bee
 \a_{ij}v_{k}+\a_{jk}v_{i}+\a_{ki}v_{j}&=&0,\  \label{idPl1}\\
 \a_{ij}+a_{ji}&=&0,\label{idPl2}\\
 \a_{ij}\a_{kl}+\a_{ik}\a_{lj}+\a_{il}\a_{jk}&=&0.\label{idPl3}
 \eee
 
 The following lemma is well known (see, for instance, \cite{Muk}). 
 \begin{lem}\label{lem1}
 The algebra $S_n$ is the free algebra modulo relations \eqref{idPl2},\eqref{idPl3}. Moreover, it has the following base over $F$:
 \bee\label{baseS}
 B_n=\{\a_{i_1j_1}\a_{i_2j_2}\cdots \a_{i_rj_r}\,|\, i_1\leq i_2\leq\cdot\leq i_r,\  j_1\leq j_2\leq\cdot\leq j_r;\, i_s<j_s\}.
 \eee
 \end{lem}
  In fact, the algebra $S_n=F[Gr(2,n)]$ is the coordinate algebra of  grassmanian $Gr(2,n)$ (see, for example, \cite[vol.1, p.42]{Sha}).
 \smallskip
 
 Let $S_{n,m}=F[\a_{ij}\,|\, j\geq m]$, then we have
 \bes 
 S_{n,n}\subseteq S_{n,n-1}\subseteq\ldots\subseteq S_{n,2}=S_{n,1}=S_n.
 \ees
 Denote by $I_m$ be the ideal of $S_n$ generated by the set $\{\a_{ij}\,|\,i<j<m\},\, 3\leq m\leq n$; then $0=I_2\subseteq I_3\subseteq I_4\subseteq \cdots\subseteq I_n$. Let also $\bar S_{n,m}=(S_{n,m}+I_m)/I_m$.
\begin{lem}\label{lem2}
The image in $\bar S_{n,m}$ of the following set forms a base of the algebra $\bar S_{n,m}$ over $F$:
\bes
B_{n,m}=\{\a_{i_1j_1}\a_{i_2j_2}\cdots \a_{i_rj_r}\in B_n\,|\, j_1\geq m\}.
\ees 
\end{lem}
{\bf Proof.}
It is easy to prove using identities \eqref{idPl3} that the set $B_{n,m}$ spans  $S_{n,m}$ modulo $I_m$.
Let us prove that it is linearly independent modulo $I_m$.   
Consider the algebra $E=E_{n,m}=F[x_m,\ldots,x_n;t;y_1,\ldots,y_n]$ and  the homomorpism 
\bes
\phi:A_n\rightarrow E,\ y_i\mapsto y_i,\, i=1,\ldots,n;\, x_j\mapsto x_j,\, j\geq m; x_k\mapsto ty_k,\,k<m.
\ees
 Let $D_{n,m}=\phi(S_n)$. Note that $\phi(\a_{ij})=0$ if $j<m$, hence $I_m\subseteq \ker\phi$. 
Introduce the $\deg{\rm lex}$ order in $E$ by setting 
\bes
x_n>x_{n-1}>\cdots>x_m>y_1>y_2>\cdots>y_n>t,
\ees
 and let $\bar f$ denotes the leading term of polynomial $f$. Then  for $i<m\leq j$ we have $\overline{\phi(\a_{ij})}=\overline{ty_iy_j-x_jy_i}=-x_jy_i$, hence
\bes
 \overline{\phi(\a_{i_1j_1}\a_{i_2j_2}\cdots \a_{i_sj_s})}=(-1)^sx_{j_1}\cdots x_{j_s}y_{i_1}\cdots y_{i_s}.
 \ees
Therefore, if $u$ and $v$ are monomials in $\phi(\a_{ij})$ then $u=v$ if and only if $\bar u=\bar v$. This easily implies that the set $\phi(B_{n,m})$ is linearly independent over $F$, and thus  the set $B_{n,m}$ is linearly independent modulo $I_m$.

\hfill$\Box$ 

Consider the elements of the base $B_n$ with more details. For any $u\in B_n$ of form \eqref{baseS} there exist uniquely defined numbers $ l<p$ such that 
\bee\label{idIS} 
u= \a_{i_1j_1}\cdots \a_{i_lj_l}\cdots \a_{i_pj_p}\cdots \a_{i_rj_r}, 
\eee
where $ j_l<m, j_{l+1}\geq m;\ i_{p}<m, i_{p+1}\geq m$.
\begin{lem}\label{lem3}
The intersection $I_m\cap S_{n,m}$ has a base formed by elements \eqref{idIS} with $r-p\geq l\geq 1$.
\end{lem}
{\bf Proof.}  Let us first prove that every element $u$ of form  \eqref{idIS} with $r-p\geq l\geq 1$ belongs to $I_m\cap S_{n,m}$. Since $l>1,\ u\in I_m$. In order to prove that $u\in S_{n,m}$, it suffices to show that
\bes
(\a_{i_1j_1}\cdots \a_{i_lj_l})( \a_{i_{p+1},j_{p+1}}\cdots \a_{i_rj_r})\in S_{n,m}.
\ees
Since $r-p\geq l$, it suffices to prove that every product $\a_{ij}\a_{kl}$ with $j<m,\,k>m$  belongs to $S_{n,m}$. But this follows easily from relation \eqref{idPl3}.

In order to prove the inverse inclusion, we associate with any element $u\in B_n$ the set of its indices ${\rm ind}(u)=\{i_1,j_1,\ldots,i_r,j_r\}$.  Note that relation \eqref{idPl3}  does not change the set of indices, hence the algebras $S_n$ and $S_{n,m}$ are homogeneous with respect to the sets of indices, that is, they may be represented as direct sums of subspaces 
with the same sets of indices.  Moreover, so is the ideal $I_m$.  Since the elements of $S_{n,m}$ are polynomials in $\a_{ij}$ with $j\geq m$, it is clear that for any homogeneous element $u\in S_{n,m}$ with ${\rm ind}(u)=\{i_1,j_1,\ldots,i_r,j_r\}$ we should have at least $r$ indices that are greater or equal to $m$.
Assume now that $\sum \l_i u_{i}\in  I_m\cap S_{n,m}$ for some $u_{i}$ of form \eqref{idIS}, then we have $\ind(u_i)=\ind(u_j)=\{i_1,j_1,\ldots,i_l,j_l,\ldots,i_p,j_p,\ldots,i_r,j_r\}$ with $l\geq 1$ for all the summonds  $u_i,u_j$. The set $\ind(u_i)$ has $l+p$ indices which are smaller than $m$, hence the sum lies in $S_{n,m}$ only if $l+p\leq r$ or $r-p\geq l$. 

\hfill$\Box$ 

 \begin{lem}\label{lem4}
  $(S_{n,m}v_m)\cap (\sum_{j<m} S_{n,j}v_j)=(S_{n,m}\cap I_m)v_m.$ 
 \end{lem}
{\bf Proof.} 
Let us first prove that  $(S_{n,m}\cap I_m)v_m\subseteq \sum_{j<m} S_{n,j}v_j$. Let $u$ be an element of form \eqref{idIS} with $r-p\geq l\geq 1$, then $u=\a_{i_1j_1}u'$, where $i_1<j_1<m$ and $u'$ is an element of form \eqref{idIS} with  $r'=r-1,\,p'=p-1,\,l'=l-1$. In particular, we have $r'-p'=r-p\geq l>l'$, therefore as in the proof of lemma \ref{lem3} we have $u'\in S_{n,m}$. Now by \eqref{idPl1} 
\bes
uv_{m}=u'(\a_{i_{1}j_{1}}v_{m})=-u'(\a_{i_1m}v_{j_1}-\a_{j_1m}v_{i_1})\in S_{n,m}v_{j_1}+S_{n,m}v_{i_1}\subseteq \sum_{j<m} S_{n,j}v_j.
\ees
Note that lemma \ref{lem2} implies that $S_{n,m}=(S_{n,m}\cap I_m)\oplus F\cdot B_{n,m}$.
Let us prove that $ B_{n,m}V_m\cap (\sum_{j<m}S_{n,j}v_j)=0$.  Assume that $w_1=w_2\neq 0$, where
\bes
w_1=\sum_{a_i\in B_{n,m}}\l_ia_iv_m,\ \ \ w_2=\sum_{j<m,\, b_j\in S_{n,j}}\mu _jb_jv_j;\  \ \l_i,\mu_j\in F.
\ees
In particular, we have $\sum_i \l_i a_ix_m=\sum_j \mu_jb_jx_j\neq 0$. Consider the leading terms of both parts with respect to the $\deg{\rm lex}$ order in $A_{n}$ when
\bes
x_n>x_{n-1}>\cdots>x_1>y_1>y_2>\cdots>y_n.
\ees
We have 
\bes
\overline{\sum_i \l_i a_ix_m}&=&\overline{\a_{i_1j_1}\ldots\a_{i_sj_s}x_m}=f(y)x_mx_{j_1}\cdots x_{j_s},\ m\leq j_1\leq j_2\leq\cdots\leq j_s,\\
\overline{\sum_j \mu_j b_jx_j}&=&\overline{\a_{p_1q_1}\ldots\a_{p_sq_s}x_j}=g(y)x_jx_{q_1}\cdots x_{q_s},\ j<m, j\leq q_1\leq q_2\leq\cdots\leq q_s.
\ees
Since $\overline{\sum_i \l_i a_ix_m}\neq \overline{\sum_j \mu_j b_jx_j}$, we conclude that $w_1=w_2=0$.

\smallskip

Now, since  $(S_{n,m}\cap I_m)v_m\subseteq \sum_{j<m} S_{n,j}v_j$, we have 
\bes
(S_{n,m}v_m)\cap(\sum_{j<m}S_{n,j}v_j)&=&((S_{n,m}\cap I_m )\oplus F\cdot B_{n,m})v_m\cap(\sum_{j<m}S_{n,j}v_j)\\
&=&((S_{n,m}\cap I_m )v_m+B_{n,m}V_m)\cap(\sum_{j<m}S_{n,j}v_j)\\
&=&(S_{n,m}\cap I_m )v_m+ B_{n,m}V_m\cap(\sum_{j<m}S_{n,j}v_j)\\
&=&(S_{n,m}\cap I_m )v_m.
\ees
\hfill$\Box$
 \begin{lem}\label{lem5}
 The subspace $S_nV\subseteq A_n^2$ is decomposed into a vector space direct sum  
 \bes
 S_nV=\oplus_{i=1}^n B_{n,i}V_i.
 \ees
 \end{lem}
{\bf Proof.}  First of all, note that due to \eqref{idPl1} we have $S_nV=\sum_{i=1}^nS_{n,i}v_i$. Denote $U_m=\sum_{i=1}^mS_{n,i}v_i$,
then $U_1\subseteq U_2\subseteq\cdots\subseteq U_n=S_nV$. Furthermore, let $W_i=U_i/U_{i-1}$, then we have a vector space isomorphism
$S_nV\cong \oplus_{i=1}^nW_i$. Finally, for any $m\leq n$ we have
\bes 
W_m&=&(\sum_{i=1}^mS_{n,i}v_i)/(\sum_{i=1}^{m-1}S_{n,i}v_i)\cong S_{n,m}v_m/(S_{n,m}v_m\cap (\sum_{i=1}^{m-1}S_{n,i}v_i))\\
&=& \hbox{ (by lemma \ref{lem4}) } = S_{n,m}v_m/(S_{n,m}\cap I_m)v_m\cong(S_{n,m}/(S_{n,m}\cap I_m))v_m\\
&\cong&  \hbox{ (by lemma \ref{lem2}) } \cong B_{n,m}V_m.
\ees
\hfill$\Box$

\begin{thm*}\label{thm7}
The space $S_n+S_nV$ is a subalgebra of algebra $\G(A_n)$ which is isomorphic to the free $\G$-algebra $\G[\emptyset;V]$.
It has a base $B_n\cup (\cup_{j=1}^n B_{n,j}v_j)$.
\end{thm*}
{\bf Proof.}
Consider the epimorphism $\f:\G[\emptyset;V]\rightarrow S_n+S_nV$ defined by the conditions $v_i\mapsto (x_i,y_i)$. Relations \eqref{idPl1} -- \eqref{idPl3} hold in the algebra $G[\emptyset;V]$ as well, and using these relations it is easy to see that it is spanned by the set $B_n\cup (\cup_{j=1}^n B_{n,j}v_j)$.
Since its image is linearly independent in $\G(A_n)$, it forms a base of $G[\emptyset;V]$, and $\f$ is an isomorphism.

\hfill$\Box$

\begin{thm*}\label{thm8}
The free $\G$-algebra $\G[t_1,\ldots,t_m;v_1,\ldots,v_n]$  on even generators $t_1,\ldots,t_m$ and odd generators $v_1,\ldots,v_n$ has the following structure:
\bes
\G_0&=&F\la t_1,\ldots,t_m\ra +S_n'\otimes F[\bar t_1,\ldots,\bar t_m],\\
 \G_1&=&F[\bar t_1,\ldots,\bar t_m]\otimes (\oplus_{j=1}^nB_{n,j}V_j),
\ees
where $F\la t_1,\ldots,t_m\ra$ and $F[\bar t_1,\ldots,\bar t_m]$ are the free associative  and the polynomial  algebras on $m$ variables, 
$S_n'$ stands for the augmentation ideal of the algebra $S_n$, $v_i\cdot v_j=\a_{ij}\in S_n,\ V_i=Fv_i$, and for any $f=f(t_1,\ldots,t_m)\in F\la t_1,\ldots,t_m\ra$ and $v\in \G_1$,   $f\cdot v=f(\bar t_1,\ldots,\bar t_m)\otimes v$.
\end{thm*}

\section{Acknowledgements}
The paper was written during I.\,Shestakov's visite to the Shenzhen  International Mathematical Center of Southern University of Science and Technology (SUSTech). He thanks Professor Efim Zelmanov for the invitation and  the Mathematical department of  SUSTech for the support and hospitality.  He was also supported by the Brazilian grants FAPESP  2018/23690-6 and CNPq 304313/2019-0.  A.Grishkov acknowledges financial support from FAPESP,  grant 2018/23690-6 (Brazil)(1-3) and from RSF,  grant 22-21-00745 (Russia)(4-5).

\bigskip

\noindent
A.\,Grishkov, University of S\~ao Paulo (S\~ao Paulo,  Brazil) and Dostoevsky Omsk State University (Omsk, Russia).\\
{\em e-mail: shuragri@gmail.com}.

\smallskip
\noindent
I.\,Shestakov,  University of S\~ao Paulo (S\~ao Paulo, Brazil),  Sobolev Institute of Mathematics (Novosibirsk, Russia), and SICM, Southern University of Science and Technology (Shenzhen, China).\\
{\em e-mail: shestak@ime.usp.br}.
\end{document}